\documentclass[hein,final,leqno]{shein}
\usepackage{latexsym}
\usepackage{graphics}
\usepackage{epsfig}
\usepackage{epsf}
\usepackage{eucal}
\usepackage{pstcol}
\usepackage{amsfonts}
\usepackage{amstext}
\usepackage{amssymb}
\usepackage{amsmath}
\usepackage{amsbsy}
\usepackage{amsxtra}
\usepackage{amscd}
\newtheorem{ltheorem}{Theorem}           

\renewcommand{\abstractname}{Abstract}
\begin{document}
\title{A DSC approach to Computational Fluid Dynamics}
\author{Steffen Hein
}                     
%
%
\institute{SPINNER GmbH., M\"unchen, Theoretical Numerics [TB01],
Aiblinger Str. 30,\linebreak DE-83620 Westerham,
Germany  \hfill \email{s.hein@spinner.de}}
%
%
\date{}
%

\maketitle 
%
\begin{abstract}
This paper presents the Dual Scattering Channel (DSC) numerical
solution of the \textsc{Navier-Stokes} equations for viscous
quasi-incompressible flow in the \textsc{Oberbeck-Boussinesq}
approximation. The implementation in hexahedral non-orthogonal
mesh is outlined. A numerical example illustrates the approach.
{\vspace{.05cm}\hfill
\textbf{MSC-classes}:\;\textnormal{65C20,\,65M06,\,76D05}}

\vspace{.2cm}
\hspace{-.52cm}
\!\!Westerham, March 31, 2006

\end{abstract}

\markboth{{\normalsize \textsc{Steffen Hein}}}
{{\normalsize \textsc{DSC approach to Fluid Dynamics}}}

\normalsize
\vspace{-0.2cm}
\section{Introduction}
\vspace{-0.2cm}
Dual Scattering Channel (DSC) schemes are characterized by a two-step
cycle of iteration which alternately updates the computed fields within
cells and on their interfaces. If the updating instructions are explicit,
then a near-field interaction principle leads to the typical structure
of the DSC algorithm and to a scattering process interpretation.
A pair of vectors that represent the same field within a cell
and on its surface essentially constitutes a \emph{scattering channel}.
Equivalently, scattering channels are sometimes defined as pairs of
distributions that 'measure' the field within the cell and on one of
its faces. A well known DSC scheme is the Transmission Line Matrix (TLM)
numerical method along \textsc{Johns}' line \cite{JoB} wherein
transmission line links visualise the scattering channels.
The TLM picture of wave propagation fails in Computational Fluid Dynamics,
thus giving rise to the far more general DSC setup.

In technical depth, DSC schemes and their relation to the TLM method
have been analysed in \cite{He1}.
Another study shows that they are unconditionally stable
under quite general circumstances, made tangible with the notion of
\emph{$\alpha$-passivity} \cite{He2}.
DSC schemes are particularly suitable for handling boundary conditions, 
non-orthogonal mesh, or also for replacing a staggered grid where
otherwise need of such is.

In section~\ref{S:sec1} we first recapitulate some characteristic
features of DSC schemes, which in extenso are treated in \cite{He1},
before we sketch in section~\ref{S:sec2} the
\textsc{Oberbeck-Boussinesq} approximation
to the \textsc{Navier-Stokes} equations.
The DSC model then outlined in section~\ref{S:sec3} represents a
\emph{prototype\/}, first of all.
In fact, the \textsc{Boussinesq} equations for viscous quasi-incompressible
flow, inspite of retaining the often predominant non-linear advective term
of the \textsc{Navier-Stokes} momentum equations,
can only claim limited range of validity,
due to their well known simplifications.
The \textsc{Oberbeck-Boussinesq} approach is yet prototypical also
in providing a basis for many turbulence models \cite{ATP}
which can be implemented following the general guide lines of this
paper.

The algorithm with non-orthogonal hexahedral mesh cell outlined
in section~\ref{S:sec3} has only recently been implemented by
\textsc{Spinner}, where it is directly coupled to a Maxwell field solver,
allowing thus for computing conductive and convective heat transfer
simultaneously with the heat sources of a lossy Maxwell field.
Besides the underlying ideas and simplifications that enter the
\textsc{Oberbeck-Boussinesq} approximation and its DSC formulation, 
a set of numerical results are displayed to illustrate and validate
the approach.
\vspace{-0.2cm}
\section{Elements of DSC schemes}\label{S:sec1}
\vspace{-.2cm}
Given a mesh cell, a \emph{port} is a vector valued
distribution associated to a cell face which assigns a state vector
$z^{\,p}\,=\,(\,p\,,\,Z\,)\,$
to a physical field $Z\,$ (any suitable smooth vector
valued function in space-time). 

We also require that, within the given cell, a
\emph{nodal image} ${p\sptilde\/}$ of $p\/$ exists,
such that
\begin{equation}\centering\label{1.1}
(\,p{\;\sptilde},\,Z\,)\;
=\;(\,p\,\circ\,\sigma\,,\,Z\,)\;
=\;(\,p\,,\,Z\,\circ\,\sigma^{-1}\,)\;,
\end{equation}
for every $\,Z\,$ (\,of class $\,C^{\,\infty}\,$, e.g.\,)\;,
where $\,\sigma\,$ denotes the spatial translation
${\sigma:\mathbb{R}^3\,\to\,\mathbb{R}^3}$ 
that shifts the geometrical node (centre of cell) onto the
(centre of) the respective face, cf. Fig1.
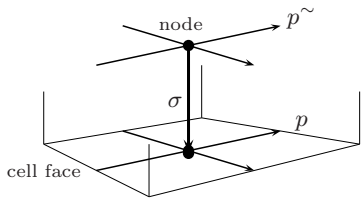
\begin{figure}[!h]\centering
\vspace{-2.3cm}
\setlength{\unitlength}{.7cm}
\begin{pspicture}(0.0,0.0)(6.0,4.6)\centering
\psset{xunit=0.7cm,yunit=0.7cm}
\psline[linewidth=0.1mm]{-}(0.0,1.0)(2.0,0.0)
\psline[linewidth=0.1mm]{-}(2.0,0.0)(6.0,1.0)
\psline[linewidth=0.1mm]{-}(6.0,1.0)(3.0,1.5)
\psline[linewidth=0.1mm]{-}(3.0,1.5)(0.0,1.0)
\psline[linewidth=0.1mm]{-}(0.0,1.0)(0.0,2.0)
\psline[linewidth=0.1mm]{-}(2.0,0.0)(2.0,1.0)
\psline[linewidth=0.1mm]{-}(6.0,1.0)(6.0,2.0)
\psline[linewidth=0.1mm]{-}(3.0,1.5)(3.0,2.5)
\psline[linewidth=0.2mm]{->}(1.0,0.5)(4.5,1.25)
\psline[linewidth=0.2mm]{->}(1.5,1.25)(4.0,0.5)
\psline[linewidth=0.2mm]{->}(1.0,2.5)(4.5,3.25)
\psline[linewidth=0.2mm]{->}(1.5,3.25)(4.0,2.5)
\psline[showpoints=true,linewidth=0.4mm]{->}(2.75,2.875)(2.75,0.825)
\psline[showpoints=true,linewidth=0.4mm]{-}(2.75,0.875)(2.75,0.825)
\rput(2.6,3.3){\scriptsize node}
\rput(4.9,3.4){\small $p\sptilde$}
\rput(2.5,1.8){\small $\sigma$}
\rput(0.0,0.5){\scriptsize cell face}
\rput(4.9,1.4){\small $p$}
\end{pspicture}
\caption{\textsl{Port on a cell face with nodal image.}\qquad\qquad\qquad\qquad
\qquad\qquad\qquad\hfill}\label{F:1}
\vspace{-0.2cm}
\end{figure}

DSC fields split thus into \emph{port} and \emph{node} components,
$z^p$ and $z^n\,$, which represent the field on the cell surfaces
and within the cells.
The two components are updated at even and odd integer multiples,
respectively, of half a timestep ~${\tau}\,$ and are usually 
constantly continued as step functions over the subsequent time
intervals of length $\tau\,$

Moreover, we assume that the updating instructions are \emph{explicit\/},
i.e., with possibly time dependent functions ${F}$ and ${G}\,$, 
for $t=m\tau\,;\,m\in\mathbb{N}\;$
\begin{equation}\centering\label{1.2}
\begin{split}
\begin{aligned}
&z^n\,(\,t\,+\,\frac{\tau}{2}\,)\;
&&:\,=\quad F\,(\,[\,z^p\,]_{\,t\,}\,,
\;[\,z^n\,]_{\,t\,-\,\frac{\tau}{2}}\,)\;,\\
&z^p\,(\,t\,+\,\tau\,)\;
&&:\,=\quad G\,(\,[\,z^p\,]_{\,t}\,,
\;[\,z^n\,]_{\,t\,+\,\frac{\tau}{2}}\,)\;,
\end{aligned}
\end{split}
\end{equation}
where ${[\,z\,]_{\,t}}$ stands for the entire sequence up to time $\,t$
\begin{equation}\label{1.3}\notag
[\,z\,]_{\,t}\,:\,=\;(\,z\,(\,t\,-\,\mu\tau\,)\,)\,_{\mu\,\in\,\mathbb{N}}\quad
\end{equation}
and we agree upon fixing
${\,z^{\,p,n}\,(\,t\,)\,:\,=\,0\,}$ for ${\,t\,<\,0\,}$.
(\,The 'back in time running' form of the sequence has certain technical 
advantages, cf.~\cite{He1}.) 

A fundamental DSC principle is \emph{near-field interaction},
which is related to the \textsc{Courant-Levi} stability criterion 
in spelling that every updated state depends only on states
(up to present time $\,t\,$) in the immediate neighbourhood.
More precisely: The next nodal state depends only on states
(along with their history) in the same cell and on its boundary,
and a subsequent port state depends only on states
(with history) on the same face and in the adjacent nodes.

As a consequence of near-field interaction, every DSC process
allows for an interpretation as a multiple scattering process
in the following sense.

If $M$ is a mesh cell system and
${\,\partial\zeta}$
denotes the boundary of cell ${\zeta\in M}$, then every DSC state
permits a unique \emph{scattering channel representation\/} in the space
\begin{equation}\centering\label{1.4}\notag
P\;:\,=\;\prod\nolimits_{\,\zeta\in M}\;
\prod\nolimits_{\,p\in\partial\zeta}\;
(\, z_{\,\zeta}^{\,p}\,,\;z_{\,\zeta}^{p\sptilde}\,)\;,
\end{equation}
with canonical projections
${\pi_{\,\zeta}^{\,p,\,n}\,:\,P\,\to\,P_{\,\zeta}^{\,p,\,n}}$
into the port and node components of cell $\zeta\,$
(\,the cell index can usually be omitted without danger of confusion\,).
Also, there is a natural involutary isomorphism 
${nb\,:\,P\,\to\,P\,}$
\begin{equation}\centering\label{1.5}\notag
nb\,:\,(\,z^{\,p}\,,\,z^{\,p\sptilde}\,)\,\mapsto\,
(\,z^{\,p\sptilde},\,z^{\,p}\,)\;,
\end{equation}
which is called the \emph{node-boundary map\/} and obviously 
maps ${\,P^{\,p}}$ onto ${\,P^{\,n}}$ and vice versa.
For every DSC process $\,z\,=\,{\,(\,z^{\,p}\,,\,z^{\,n}\,)(\,t\,)\,}$,
the following \emph{incident\/} and \emph{outgoing fields}
$\,{z_{\,in}^{\,p}}\,$ and $\,{z_{\,out}^{\,n}}\,$ 
are then recursively well defined, and are processes in
${\,P^{\,p}\,}$ and ${\,P^{\,n}\,}$, respectively:\; \\ 
For $\,t < 0\,$, $\,{z_{\,in}^{\,p}\,(\,t\,)\,
:\,=\,z_{\,out}^{\,n}\,(\,t-\frac{\tau}{2}\,)\,:\,=\,0\;}\,$,
and for every $\,0\,\leq\,t\,=\,m\tau\,$;
$\,m\in\mathbb{N}\,$
\begin{equation}\centering\label{1.6}
\begin{split}
\begin{aligned}
z_{\,in}^{\,p}\,(\,t\,)\;
&:\,=\; z^{\,p}\,(\,t\,)\,
-\,nb\circ z_{\,out}^{\,n}\,(\,t\,-\,\frac{\tau}{2}\,)\;,\\
z_{\,out}^{\,n}\,(\,t\,+\,\frac{\tau}{2}\,)\;
&:\,=\; z^{\,n}\,(\,t+\,\frac{\tau}{2}\,)\,-\,nb\circ z_{\,in}^{\,p}\,(\,t\,)\;.
\end{aligned}
\end{split}
\end{equation}
At every instant holds, hence, $\,{z^{\,p}\,(\,t\,)}\,
=\,{nb\circ z_{\,out}^{\,n}\,(\,t\,-\,\frac{\tau}{2}\,)}\,
+\,{z_{\,in}^{\,p}\,(\,t\,)}\,$
and \\
$\,{z^{\,n}\,(\,t\,+\,\frac{\tau}{2}\,)}\,
=\,{nb\circ z_{\,in}^{\,p}\,(\,t\,)}\,
+\,{z_{\,out}^{\,n}\,(\,t\,+\,\frac{\tau}{2}\,)}\,$\,.
Then, near-field interaction implies that every state is only a function of
states incident (up to present time $\,t\,$) on scattering channels
connected to the respective node or face.
\newline
More precisely, by induction holds
\begin{ltheorem}{\label{T:1}}.
A pair of functions
$\,\mathcal{R}\,$ and $\,\mathcal{C}\,$ exists, such that for every cell
$\,{\zeta\in M}\,$
the process
$\,z_{\,\zeta}^{\,n}\,=\,{\pi_{\,\zeta}^{\,n}\circ z}\,$
complies with
\begin{equation}\centering\label{1.9}
z_{\,\zeta}^{\,n}\,(\,t\,+\,\frac{\tau}{2}\,)\;
=\;\mathcal{R}\,(\,(\,z_{\,in}^{\,p}\,(\,t\,-\,\mu\tau\,)\,
)_{\,p\in\partial\zeta\,;\;\mu\in\mathbb{N}}\,)\;
\end{equation}
and the port process
$\,z_{\,\zeta}^{\,p}\,=\,\pi_{\,\zeta}^{\,p}\,\circ\,z\,$
satisfies
\begin{equation}\centering\label{1.10}
\begin{aligned}
z_{\,\zeta}^{\,p}\,(\,t\,+\,\tau\,)\;
=\;\mathcal{C}\,(\,(\,z_{\,out}^{\,n}\,(\,t\,+\,\frac{\tau}{2}-
\mu\tau\,)\,)_{\,n\,\mid\,\partial\zeta\,;\;\mu\in\mathbb{N}}\,)\;.\\
\scriptsize{(\;\text{'$\;\mid\,$' short-hand for 'adjacent to'}\;\;)}
\end{aligned}
\end{equation}
\end{ltheorem}
\vspace{-0.30cm}
\emph{Remarks}
\hspace{0.65cm}
\begin{itemize}
\item[(i)]
The statements immediately imply that
$\,z_{\zeta,\,out}^{\,n}\,$ and $\,z_{\zeta,\,in}^{\,p}\,$
are themselves functions of states incident
on connected scattering channels, since
\vspace{-.3cm}
\end{itemize}
\begin{equation}\centering\label{1.11}
\begin{aligned}
z_{\zeta,\,out}^{\,n}\,(\,t\,+\,\frac{\tau}{2}\,)\;
=\;\mathcal{R}\,(\,(\,z_{\,in}^{\,p}\,(\,t\,-\mu\tau\,)\,)_{\,p\,\in\,
\partial\zeta\,;\;\mu\in\mathbb{N}}\,)\,
-\,z_{\zeta,\,in}^{\,p}\,(\,t\,) \quad \text{and} \\
z_{\zeta,\,in}^{\,p}\,(\,t\,)\;
=\;\mathcal{C}\,(\,(\,z_{\,out}^{\,n}\,(\,t\,-\,\frac{\tau}{2}-\mu\tau\,)\,
)_{\,n\,\mid\,p\,;\;\mu\in\mathbb{N}}\,)\,
-\,z_{\zeta,\,out}^{\,p\sptilde}\,(\,t\,-\,\frac{\tau}{2}\,)
\end{aligned}
\end{equation}
\begin{itemize}
\item[(ii)]
$\,\mathcal{R}\,$ and $\,\mathcal{C}\,$ are named the
\emph{reflection} and \emph{connection\/} maps, respectively,
of the DSC algorithm.
\end{itemize}
\begin{itemize}
\item[(iii)]
Near field interaction implies computational stability, if the reflection
and connection maps are contractive or 
\nolinebreak{$\,\alpha$-\emph{passive}} \cite{He2},
in addition.
\end{itemize}
\vspace{-.2cm}
\vspace{-0.2cm}
\section{The dynamic equations}\label{S:sec2}
\vspace{-0.2cm}
The DSC algorithm is thus simply characterized as a two-step
explicit scheme that alternately updates states in ports and nodes
of a cellular mesh and which, in virtue of a near-field interaction
principle, allows for a canonical interpretation as a scattering process,
the latter exchanging \emph{incident} and \emph{reflected} quantities
between cells and their interfaces.

Ports and nodes are related to physical fields by vector valued
distributions that evaluate the fields on cell faces and within
cells of a cellular mesh.
Such a distribution may be a \emph{finite integral}
-~as is the case with the TLM method, where finite path integrals over
electric and magnetic fields are evaluated in a discrete approximation
to Maxwell's integral equations \cite{He3}.
In the simplest case, it is a \emph{Dirac measure} that pointwise
evaluates a field (\,or a field component\,) within a cell and on its
surface. The distribution can also be a composite of Dirac measures
that evaluate a field at different points in the cell -~which applies,
for instance, to the gradient functional treated in sect.~\ref{S:sec3}.

Classical thermodynamics with, in particular, energy conservation
entail the \emph{convection-diffusion} equation for the temperature
$\,T\,$ in a fluid of velocity $\,\vec{u}\,$ with constant thermal
diffusivity $\,\alpha\,$, heat source(s) $\,q\,$,
and negligible viscous heat dissipation, viz.
\vspace{-0.15cm}
\begin{equation}\centering\label{2.1}
\frac{\partial\,T}{\partial\,t}\;
+\;\vec{u}\,\cdot\,grad\;T\;
=\;\alpha\,\Delta\,T\;+\;q\;.
\vspace{-0.10cm}
\end{equation}
This is the energy equation for \emph{Boussinesq-incompressible\/}
fluids, e.g. \cite{GDN}.
The Navier-Stokes \emph{momentum equations\/} for a fluid of
dynamic viscosity $\,\mu\,$,
under pressure $\,p\,$,
and in a gravitational field of acceleration $\,\vec{g}\,$ require
\vspace{-0.15cm}
\begin{equation}\centering\label{2.2}
\frac{\partial}{\partial\,t}\,(\,\varrho\,\vec{u}\,)\;
+\;(\,\vec{u}\cdot grad\,)\,(\,\varrho\,\vec{u}\,)
+\;grad\,p\;
=\;\mu\,\Delta\,\vec{u}\;
+\;\varrho\;\vec{g}\;.
\vspace{-0.10cm}
\end{equation}
The \emph{Oberbeck-Boussinesq\/} approximation \cite{Obb},\cite{Bss}
starts from the assumption that the fluid properties are constant,
except fluid density, which only in the gravitational term varies
linearly with temperature;
and that viscous dissipation can be neglected.
Equations \eqref{2.2} become so with
$\,\varrho_{\,\infty}=const\,$ and
$\,\varrho\,(\,T\,)=\varrho_{\,\infty}\,\beta\,
(\,T\,(\,t\,,\vec{x}\,)\,-\,T_{\,\infty}\,)\;$; $\,\beta\,
:\,=\,\varrho^{-1}\,\partial\,\varrho\,/\,\partial\,T\,$
\vspace{-0.15cm}
\begin{equation}\centering\label{2.3}
\frac{\partial\,\vec{u}}{\partial\,t}\,
+\,(\vec{u}\cdot grad\,)\,\vec{u}\,
+\,\frac{grad\,p}{\varrho_{\infty}}\,
=\,\frac{\mu}{\varrho_{\infty}}\,\Delta\,\vec{u}\;
+\,\beta\,(\,T\,(\,t,\,\vec{x}\,)
- T_{\infty}\,)\,\vec{g}\,.
\vspace{-0.10cm}
\end{equation}
Applying the Gauss-Ostrogradski theorem to the
integrals over $\,{\Delta}\,=\,{div\,grad\,}\,$
on cell $\zeta$ with boundary ${\partial\zeta}\/$
and using 
$\vec{u}\,\cdot\,grad\,f\,+\,f\,div\,\vec{u}\,=
\,div\,(\,f\,\vec{u}\,)\,$,
equations (\ref{2.1}, \ref{2.3}) yield,
with a time increment $\tau\,$,
the following updating instructions for nodal $T\,$ and $\vec{u}$,
the latter, along with $q$\,, averaged over the cell volume 
${V_{\zeta}}$
\vspace{-0.15cm}
\begin{equation}\centering\label{2.4}
\begin{aligned}
T\,(\,t+\frac{\tau}{2}\,)\,
&:\,=\,T\,+\,\tau\,(\;T\,div\;\vec{u}\,+\,q\,\;)\;+\\
&+\;\frac{\tau}{V_{\zeta}}\,\int\nolimits_{\partial\,\zeta}(\,\alpha\,grad\,T\,
-\,T\,\vec{u}\,)\cdot\,dF\,
\end{aligned}
\vspace{-0.3cm}
\end{equation}
and
\vspace{-0.3cm}
\begin{equation}\centering\label{2.5}
\begin{aligned}
\vec{u}\,(\,t+\frac{\tau}{2}\,)\,
&:\,=\,\vec{u}\,+\,\tau\,(\;\vec{u}\,div\,\vec{u}\,
-\,\beta\,(\,T\,-\,T_{\infty}\,)\,\vec{g}\,
-\,\frac{grad\,p}{\varrho_{\infty}}\,\;)\;+\\
&+\;\frac{\tau}{V_{\zeta}}\,\{\;
\int\nolimits_{\partial\,\zeta}\frac{\mu}{\varrho_{\infty}}\
grad\,\vec{u}\cdot\,dF\,
-\,\int\nolimits_{\partial\,\zeta}\vec{u}\;(\,\vec{u}\cdot dF\,)\;\}
\end{aligned}
\vspace{-0.10cm}
\end{equation}
In equations \eqref{2.4} and \eqref{2.5} the last former updates 
(at time $t-\tau/2\/$) of the nodal quantities enter the right-hand
sides in the first line, and the former updates (at time $t\/$)
of the cell face quantities enter the second line.
Note that $div\,\vec{u}\/$ vanishes, of course, for incompressible flow.

All nodal values of $\,T\,$ and $\,\vec{u}\,$ are thus updated
at the reflection step of the algorithm, while the cell face
quantities that enter the integrals at the right-hand sides
are updated on the connection step.
In next section we will see how to proceed with an unstructured 
(non-orthogonal) hexahedral mesh.
\vspace{-0.2cm}
\section{Non-orthogonal hexahedral cell}\label{S:sec3}
\vspace{-0.2cm}
The physical interpretation of a DSC algorithm associates a
smoothly varying ( e.g. in time and space $\,C^{\,\infty}$-\,)
scalar or vector field $\,Z\,$ to port and node states
$\,z^{\,p}\,$ and $\,z^{\,n}\,$ of a mesh cell system.

\begin{figure}[!h]\centering
\vspace{-0.8cm}
\setlength{\unitlength}{1.cm}
\begin{pspicture}(-1.0,-.5)(22,3.5)\centering
\psset{xunit=.4cm,yunit=.4cm}
\psline[linewidth=0.1mm]{->}(1.5,1.5)(-1.0,0.0)
\psline[linewidth=0.1mm]{->}(1.5,1.5)(5.0,2.0)
\psline[linewidth=0.1mm]{->}(1.5,1.5)(2.0,5.0)
\psline[linewidth=0.4mm]{->}(-1.0,0.0)(-1.0,4.0)
\psline[linewidth=0.4mm]{->}(-1.0,0.0)(4.0,0.0)
\psline[linewidth=0.4mm]{->}(5.0,2.0)(5.0,5.0)
\psline[linewidth=0.4mm]{->}(5.0,2.0)(4.0,0.0)
\psline[linewidth=0.4mm]{->}(4.0,0.0)(2.0,3.0)
\psline[linewidth=0.4mm]{->}(2.0,5.0)(-1.0,4.0)
\psline[linewidth=0.4mm]{->}(2.0,5.0)(5.0,5.0)
\psline[linewidth=0.4mm]{->}(-1.0,4.0)(2.0,3.0)
\psline[linewidth=0.4mm]{->}(5.0,5.0)(2.0,3.0)
\rput(0.8,0.8){$_{_{0}}e$}
\rput(5.0,1.0){$_{_{1}}e$}
\rput(3.0,4.2){$_{_{2}}e$}
\rput(0.5,4.9){$_{_{3}}e$}
\rput(4.0,1.5){$_{_{4}}e$}
\rput(3.4,5.5){$_{_{5}}e$}
\rput(0.5,2.9){$_{_{6}}e$}
\rput(1.5,-.5){$_{_{7}}e$}
\rput(1.4,4.1){$_{_{8}}e$}
\rput(-1.6,2.0){$_{_{9}}e$}
\rput(3.2,2.6){$_{_{10}}e$}
\rput(5.7,3.4){$_{_{11}}e$}
\rput(6.0,-0.5){\small{\textnormal{(a)}}}
\psline[linewidth=0.1mm]{->}(10.5,1.5)(08.0,0.0)
\psline[linewidth=0.1mm]{->}(10.5,1.5)(14.0,2.0)
\psline[linewidth=0.1mm]{->}(10.5,1.5)(11.0,5.0)
\psline[linewidth=0.1mm]{->}(08.0,0.0)(08.0,4.0)
\psline[linewidth=0.1mm]{->}(08.0,0.0)(13.0,0.0)
\psline[linewidth=0.1mm]{->}(14.0,2.0)(14.0,5.0)
\psline[linewidth=0.1mm]{->}(14.0,2.0)(13.0,0.0)
\psline[linewidth=0.1mm]{->}(13.0,0.0)(11.0,3.0)
\psline[linewidth=0.1mm]{->}(11.0,5.0)(08.0,4.0)
\psline[linewidth=0.1mm]{->}(11.0,5.0)(14.0,5.0)
\psline[linewidth=0.1mm]{->}(08.0,4.0)(11.0,3.0)
\psline[linewidth=0.1mm]{->}(14.0,5.0)(11.0,3.0)
\psline[showpoints=true,linewidth=0.4mm]{->}(12.375,3.375)(10.0,1.75)
\psline[showpoints=true,linewidth=0.4mm]{->}(09.375,2.625)(13.0,2.50)
\psline[showpoints=true,linewidth=0.4mm]{->}(11.37,.875)(11.0,4.25)
\psline[showpoints=true,
linewidth=0.6mm]{-}(11.1875,2.5625)(11.1875,2.5625) 
\rput(09.55,1.55){$_{_{0}}b$}
\rput(13.5,2.60){$_{_{1}}b$}
\rput(11.6,4.40){$_{_{2}}b$}
\rput(15,-0.5){\small{\textnormal{(b)}}}
\psline[linewidth=0.1mm]{->}(18.0,4.0)(18.0,1.0)
\psline[linewidth=0.1mm]{->}(18.0,1.0)(23.0,0.0)
\psline[linewidth=0.1mm]{->}(18.0,4.0)(22.0,5.0)
\psline[linewidth=0.1mm]{->}(22.0,5.0)(23.0,0.0)
\psline[showpoints=true,linewidth=0.5mm]{-}(20.25,2.5)(20.25,2.5) 
\psline[linewidth=0.25mm]{->}(22.5,2.5)(18.0,2.5)
\psline[linewidth=0.25mm]{->}(20.0,4.5)(20.5,0.5)
\psline[showpoints=true, linewidth=0.4mm]{->}(22.5,2.5)(26.751,3.400)
\psline[showpoints=true,linewidth=0.4mm]{->}(18.0,2.5)(16.5,2.5)
\psline[showpoints=true,linewidth=0.4mm]{->}(20.0,4.5)(19.287,7.352)
\psline[showpoints=true,linewidth=0.4mm]{->}(20.5,0.5)(19.650,-3.751)
\rput(19.0,3.2){$_{_{0}}b$}
\rput(21.0,1.5){$_{_{1}}b$}
\rput(24.6,3.6){$_{_{0}}f$}
\rput(16.8,3.2){$_{_{1}}f$}
\rput(20.2,6.4){$_{_{2}}f$}
\rput(20.7,-1.4){$_{_{3}}f$}
\rput(23.6,2.0){$_{_{0}}Z^{\,p}$}
\rput(17.1,1.8){$_{_{1}}Z^{\,p}$}
\rput(18.8,4.7){$_{_{2}}Z^{\,p}$}
\rput(19.3,0.0){$_{_{3}}Z^{\,p}$}
\rput(21.1,3.2){$Z^{\,n}$}
\rput(25,-0.5){\small{\textnormal{(c)}}}
\end{pspicture}
\caption{\textsl{Non-orthogonal hexahedral mesh cell. \newline
\textnormal{(a)} Edge vectors.\;
\textnormal{(b)} Node vectors.\;
\textnormal{(c)} Face vectors.}\hfill}\label{F:2}
\end{figure}
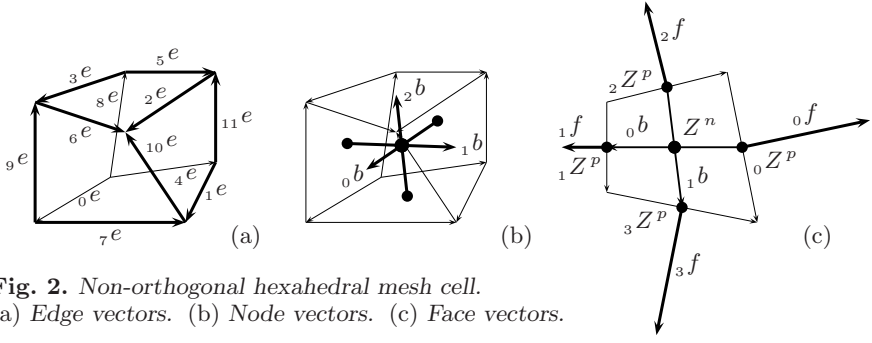

Let any hexahedral cell be given by its eight vertices.
Define then \emph{edge vectors} ${(_{\nu} e)_{\nu=0,...,11}}\/$,
\emph{node vectors\/} ${(_{\mu} b)_{\mu = 0,1,2}}$,
and \emph{face vectors\/} ${(_{\iota} f)_{\iota = 0,...,5}}$,\,
using the labelling scheme of figure~\ref{F:2}\,a
\vspace{-0.15cm}
\begin{equation}\label{3.1}\centering
\begin{split}
\begin{aligned}
_{\mu}b\;&:\,=\quad\frac{1}{4}
&& \!\!\sum\nolimits_{\nu = 0}^{3}\,_{_{(4\mu+\nu)}}e\,\,
&&\mu\,=\,0,1,2\\
\text{and}\quad
_{\iota}f\;&:\,=\;\,\frac{(-1)^{\,\iota}}{4}
&&\,(\,\,_{_{(8+2\iota)}}e\,
+_{_{(9+2(\iota+(-1)^{\iota}))}}e\,)
\,\,\land &&\\
& &&\;\;\land\,(\,_{_{(4+2\iota)}}e\, 
+_{_{(5+2\iota)}}e\,)\,\,&&\iota\,=\,0,...,5\,,
\end{aligned}
\end{split}
\vspace{-0.10cm}
\end{equation}
with all indices understood cyclic modulo 12\, and $\,\land\,$
denoting the cross product in $\mathbb{R}^3$.

At every cell face ${\iota\in\{0,...,5\}}\/$ and for any 
given $\tau\in\mathbb{R}_{+}\,$ the following time shifted finite
differences of $\,Z\,$ in directions ${ _{\mu} b }$ (\,$\mu = 0,1,2\,$)
form a vector valued function
\vspace{-0.15cm}
\begin{equation}\label{3.2}\centering
\begin{split}
_{\iota}\!{\nabla}^{B} Z_{\mu}\,(\,t\,)\;:\,=\;
\begin{cases}
\,2\,(-1)^{\iota}(\,Z^{n}\, _{\mid\,t-\tau/2} 
-\,_{\iota}Z^{p}\,_{\mid\,t}\,)\quad
&\text{if $\mu\,=\,[\iota / 2]$}\\
\,(\,\,_{2\mu+1} Z^{p}\,-\,_{2\mu} Z^{p}\,\,)
\,_{\mid\,t -\tau\,}\quad
&\text{if $\mu\,\neq\,[\iota / 2]$}\,
\end{cases}
\end{split}
\vspace{-0.10cm}
\end{equation}
($\,[\,x\,]$ denotes the \emph{integer part} of $\,x\in\mathbb{R}\,$).
The time increments are chosen conform with the updating conventions
of DSC schemes (as will be seen in a moment) and are consistent.
In fact, in the first order of the time increment ${\,\tau\,}$
and of the linear cell extension,
the vector ${\,_{\iota}\!{\nabla}^{B} Z\,}$
in the centre point of face ${\,\iota\,}$ approximates
the scalar products of the node vectors with the gradient
${\,\nabla Z\,}$.
Let, precisely, for a fixed centre point on face ${\,\iota\,}$
and $\,\epsilon \in \mathbb{R}_{+}\,$ the \emph{$\epsilon$-scaled cell}
have edge vectors
$\,_{\iota}e\sptilde\,:\,=\,\epsilon\,\,_{\iota}e\,$. 
Let also $\,_{\iota}{\nabla}^{B\sptilde}Z_{\mu}\,$ denote
function \eqref{3.2} for the $\epsilon$-scaled cell (with node vectors 
$\,_{\mu}b\sptilde\,=\,\epsilon\,_{\mu}b\;$).
Then at the fixed point holds
\vspace{-0.15cm}
\begin{equation}\label{3.3}\centering
\begin{split}
<\,_{\mu}b\,,\,\text{grad($Z$)}\,>\,\,\,
=\,\,_{\mu}b\cdot\nabla Z\,\,
=\,\,\lim_{\epsilon\to 0}\,\,\lim_{\tau \to 0} \, \,
\frac{1}{\epsilon}\,_{\iota}\!{\nabla}^{B^{\sptilde}}Z_{\mu} \, ,
\end{split}
\vspace{-0.10cm}
\end{equation}
as immediately follows from the required $C^1$-smoothness of the field $Z$.

To recover the gradient
${{\nabla}Z\,}$ from \eqref{3.2}
in the same order of approximation,
observe that for every orthonormal
basis $\,{(_{\nu}u)_{\nu=0,...,m-1}}\,$ of
$\mathbb{R}^{m}\,\text{or}\,\,\mathbb{C}^{m}\,$, and for
any basis $\,{(_{\mu}b)_{\mu = 0,...,m-1}}\,$ with coordinate
matrix ${\beta_{\nu}^{\mu}}\,=\,{<\, _{\nu}u\,,\, _{\mu}b\,>}$,
the scalar products of every vector $\,a\,$ with $\,{_{\mu} b}\,$ equal
\vspace{-0.15cm}
\begin{equation}\label{3.4}\centering
\underbrace{<\,_{\mu}b\,,\,a\,>}_{\qquad=\,: 
\,\,{\alpha}_{\mu}^{B}}\,\,=\,\sum\nolimits_{\nu=0}^{m-1}\,
\underbrace{<\,_{\mu}b\,,\,_{\nu}u\,>}_{\,\,\,
({\bar{\beta}}_{\mu}^{\nu})\,=\,({\beta}_{\nu}^{\mu} )^{^{*}}}\,
\underbrace{<\,_{\nu}u\,,\,a\,>}_{\qquad=\,:\,\,{\alpha}_{\nu}}
\,\,=\,\bar{\beta}_{\mu}^{\nu}\,{\alpha}_{\nu}\;
\vspace{-0.10cm}
\end{equation}
(\,at the right-hand side, and henceforth, we observe
\textsc{Einstein}'s summation convention -~without yet summing up  
over indices that anywhere appear as left-hand subscripts\,),
hence
\vspace{-0.15cm}
\begin{equation}\label{3.5}\centering
{\alpha}_{\nu}\,=\,
{\gamma}_{\nu}^{\mu}\alpha_{\mu}^{B}\,,
\qquad\text{with}\qquad({\gamma}_{\nu}^{\mu}) 
\,:\,=\,{({(\beta_{\nu}^{\mu})}^{*})}^{-1}\quad .
\vspace{-0.10cm}
\end{equation}
I.e., the scalar products of any vector with the basis
vectors ${_{\mu} b\,}$ transform into the coordinates of that vector
with respect to an orthonormal basis $\,{_{\nu}u}\,$ by multiplication
with matrix $\,{\gamma=(\beta^*)^{-1}}\,$, where
$\,{\beta_{\nu}^{\mu}}\,=\,{<\,_{\nu}u\,,\,_{\mu}b\,>}\,$,\,
i.e. $\,\beta\,$ is the matrix of the coordinate (column) vectors
$\,{_{\mu}b}\,$ with respect to the given ON-basis
$\,{_{\nu}u}\,$,\, and $\,\gamma\,$ its adjoint inverse.

This applied to the node vector basis ${ _{\mu}b\,}$ and \eqref{3.3}
yields the approximate gradient of $\,Z\,$ at face $\iota$
\vspace{-0.15cm}
\begin{equation}\label{3.6}\centering
_{\iota}\!\nabla Z_{\nu}\quad
=\quad{\gamma}_{\nu}^{\mu}\,\,\,_{\iota}\!{\nabla}^{B}Z_{\mu}.
\vspace{-0.10cm}
\end{equation}
The scalar product of the gradient with face vector  
${_{\iota}f^{\nu}}\,=\,{<\,_{\iota}f,\,_{\nu}u>}\,$,\\
$\,\nu\in\{0,1,2\}\,$ is thus
\vspace{-0.15cm}
\begin{equation}\label{3.7}\centering
\begin{split}
\begin{aligned}
\qquad _{\iota}S\; 
&=\;_{\iota}f\,\cdot\,_{\iota}\!\nabla Z\;
=\underbrace{_{\iota}f^{\nu}\;
{\gamma}_{\nu}^{\mu}}_{\qquad\;=\,:\,\,_{\iota} s^{\mu}}\,
_{\iota}\!{\nabla}^{B} Z_{\mu}\,\,
=\,\,_{\iota} s^{\mu}\,\,_{\iota}\!{\nabla}^{B} Z_{\mu}\;.
\end{aligned}
\end{split}
\vspace{-0.10cm}
\end{equation}
Continuity of the gradient at cell interfaces yields
linear updating equations for $Z^p$ on the two adjacent faces.  
In fact, for any two neighbouring cells
$\zeta$, $\chi$ with common face, labelled $\iota$
in cell $\zeta$ and $\kappa$ in $\chi\,$, continuity requires
\vspace{-0.15cm}
\begin{equation}\label{3.8}\centering
_{\iota}^{^{\zeta}}\!S\quad=\quad-\,\,_{\kappa}^{^{\chi}}\!S\,.
\vspace{-0.10cm}
\end{equation}
Substituting \eqref{3.7} for
$\,_{\iota}^{^{\zeta}}\!S\,$ and $\,_{\kappa}^{^{\chi}}\!S\,$
and observing the time shifts in \eqref{3.2}
\linebreak
provides the updating relations for 
$\,Z^{\,p}\,$ at the cell interfaces.
To make these explicit, we first introduce the quantities
$\,{_{\iota}z_{\mu}^{\,p,\,n}}\,$, ($\,\iota\,=\,0,...,5\,$;
$\mu\,=\,0,1,2\,$)
\vspace{-0.15cm}
\begin{equation}\label{3.10}\centering
\begin{split}
_{\iota}z_{\mu}^{n}\,(\,t\,)\quad :\,=\quad
\begin{cases}
\,\,2\,(-1)^{\iota}\,\,Z^{\,n}\,_{\mid\,t}\qquad
&\text{if $\mu\,=\,[\iota /2]$}\,\,\\
\,\,(\,_{2\mu +1} Z^{\,p} 
-\,_{2\mu} Z^{\,p}\,)_{\mid\,t-\tau/2}\qquad
&\text{else}\,
\end{cases}\;,
\end{split}
\vspace{-0.10cm}
\end{equation}
which in virtue of \eqref{1.1} yields
$\;{_{\iota}z_{\mu}^{\,p}}\,=\,{(\,p\,,\,Z\,)}\,
=\,{(\,p\sptilde,\,Z\circ\,_{\iota}{\sigma}^{-1}\,)}\,
=\\=\,{_{\iota}z_{\mu}^{\,n}\,\mid{Z\,\circ\,_{\iota}{\sigma}^{-1}}}\,$,
where $\,_{\iota}{\sigma}\,:\,n\,\mapsto\,p\,$
denotes the nodal shift pertinent to face $\,\iota\,$.\;
In particular
\begin{equation}\label{3.11}
_{\iota}z_{[\iota/2]}^{\,p}\,(\,t\,)\quad=\quad
\,\,2\,(-1)^{\iota}\,\,_{\iota}Z^{\,p}\,_{\mid\,t}\;,
\vspace{-0.10cm}
\end{equation}
which together with \eqref{3.10} is consistent for 
$\,\mu\,\neq\,[\iota/2]\,$ with
\begin{equation}\label{3.12}
_{\iota}z_{\mu}^{\,n}\,(\,t+\tau/2\,)\quad
=\quad -\;\frac{1}{2}\,(\,_{\,2\mu+1}z_{\mu}^{\,p}\,
+\, _{\,2\mu}z_{\mu}^{\,p}\,)\,(\,t\,)\;.
\vspace{-0.10cm}
\end{equation}
From (\,\ref{3.2}, \ref{3.7}, \ref{3.10}, \ref{3.11}\,) follows that
\vspace{-0.15cm}
\begin{equation}\label{3.13}
\begin{split}
_{\iota}S\,_{\mid\,t+\tau}\quad
&=\quad\,_{\iota}s^{\mu}\,(\,_{\iota}z_{\mu}^{\,n}\,_{\mid\,t+\tau/2}\,
-\,2\,{(-1)}^{\iota} {\delta}_{\mu}^{[\iota/2]}\,\,
_{\iota}Z^{\,p}\,_{\mid\,t+\tau}\,) \\
&=\quad _{\iota}s^{\mu}\,(\,_{\iota}z_{\mu}^{n}\,_{\mid\,t+\tau/2} \,
-\,{\delta}_{\mu}^{[\iota/2]}\,\, 
_{\iota}z_{\mu}^{\,p}\,_{\mid\,t+\tau}\,)
\; .
\end{split}
\vspace{-0.10cm}
\end{equation}
The continuity of $\,Z\,$, i.e.
$\,_{\iota}^{^{\zeta}}Z\,^{p}\,=\,_{\kappa}^{^{\chi}}Z\,^{p}\,$, 
the implies with (\ref{3.8},\ref{3.10})
\vspace{-0.15cm}
\begin{equation}\label{3.14}
_{\iota}^{^{\zeta}}z\,_{[\iota/2]}^{p}\,(\,t+\tau\,)\,=\,
\,\frac{\,_{\iota}^{^{\zeta}}s\,^{\mu}\,\, 
_{\iota}^{^{\zeta}}z\,_{\mu}^{n}\,(\,t+\tau/2\,)\,
+\,_{\kappa}^{^{\chi}}s\,^{\nu}\,\,\,
_{\kappa}^{^{\chi}}z\,_{\nu}^{n}\,(\,t+\tau/2\,)} 
{_{\iota}^{^{\zeta}}s\,^{[\iota/2]}
+\,( -1 )^{\iota +\kappa}\;_{\kappa}^{^{\chi}}s\,^{[\kappa/2]}}\;.
\vspace{-0.10cm}
\end{equation}
For completeness we agree upon setting
$\,_{\iota}^{^{\zeta}}z\,_{\mu}^{p}\,(\,t+\tau\,)\,
:\,=\,_{\iota}^{^{\zeta}}z\,_{\mu}^{n}\,(\,t+\tau/2\,)\,$
for $\,\mu\,\neq\,[\iota/2]\,$
(which yet contains a slight inconsistency,
in that continuity might be infringed; this can be remedied
by taking the arithmetic means of the two adjacent values).
 -~In fact, our agreement doesn't do harm, since any discontinuity
disappears with mesh refinement. 

We have, hence, a complete set of recurrence relations for
$\,{z^{\,p}}\,$ (\,given $\,{z^{\,n}}\,$ by the former reflection step)
which at the same time determine
the field components on face~$\iota$ and their gradients
\vspace{-0.15cm}
\begin{equation}\label{3.15}\centering
_{\iota}\!\nabla Z_{\nu}\quad
=\quad{\gamma}_{\nu}^{\mu}\; _{\iota}z_{\mu}^{\,p}.
\vspace{-0.10cm}
\end{equation}
Essentially this constitutes the connection step of the algorithm.
\newline
Nodal gradients are similarly (yet even more simply) derived using
\vspace{-0.15cm}
\begin{equation}\label{3.16}\centering\notag
{\nabla}^{B} Z_{\mu}^{\,n}\,(\,t\,+\frac{\tau}{2}\,)\quad:\,
=\quad(\,_{2\mu+1} Z^{\,p}\,-\,_{2\mu} Z^{\,p}\,)\,(\,t\,)\;;
\quad\mu\,=\,0,\,1,\,2\;
\vspace{-0.10cm}
\end{equation}
in the place of \eqref{3.2} and then again \eqref{3.6}.
With the node and cell-boundary values and gradients of $\,T\,$ and
$\,\vec{u}\,$ the nodal updating relations for the latter are
immediately extracted from equations~(\,\ref{2.4}, \ref{2.5}\,)
in sect.~\ref{S:sec2}\,. For equation \eqref{2.4} this is
essentially (up to the convective term) carried out in
\cite{He1}, sect.~5, and the procedure remains straightforward
in the case at hand.
Note that a well-timed LES coarsening routine~\cite{He3}
should be periodically carried out before the nodal step of iteration
in order to avoid instablities from the energy cascade~\cite{Po}. 

The updating relations thus obtained are explicit and consistent
with near-field interaction (\,only adjacent quantities enter\,).
So, they can optionally be transformed into scattering relations
for incident and reflected quantities \eqref{1.6} along the
guidelines of section~\ref{S:sec1} -~with established advantages
for stability estimates \cite{He2}.
\vspace{-0.2cm}
\section{Pressure}\label{S:sec4}
\vspace{-0.2cm}
Conservation of mass simply requires divergence-free flow,
$\;{div\,\vec{u}\,=\,0}\;$,
for a Boussinesq-incompressible fluid.
In integral form, using Gauss' Theorem, this means
$0\,=\,\int\nolimits_{\,\zeta}\,div\;\vec{u}\;dV\,
=\,\int\nolimits_{\,\partial\zeta}\,\vec{u}\,\cdot\,dF\,$.
Since equations (\,\ref{2.1}, \ref{2.2}\,) in section~\ref{S:sec2}
do not a priori guarantee this, additional arrangements must be
made. This is done by the following procedure, which in Magneto-Hydrodynamics
is known as \emph{divergence cleaning}.

In a successive overrelaxation (SOR) routine, carried out between
the connection and reflection steps of the iteration cycle,
firstly the (discrete) cell boundary integrals
$\,I_{\partial\zeta}\,=\,\int _{\partial\zeta}\vec{u}\cdot dF\,$
are computed and then the pressure $\,p\,$
which compensates $\,I_{\partial\zeta}\,$ so that
\vspace{-0.15cm}
\begin{equation}\centering\label{4.2}
\quad\frac{\tau}{\varrho_{\infty}}\,
\int\nolimits_{\,\partial\zeta}\, grad\,p\,\cdot\,dF\quad
=\quad\int\nolimits_{\,\partial\zeta}\,\vec{u}\,\cdot\,dF\;.
\vspace{-0.10cm}
\end{equation}
Taking indeed $\,Z\,$ of the preceeding section as the pressure,
equations \eqref{4.2}
(\,of course, in discrete form with sums over the cell faces\,)
yield a unique solution $\,p^{\,n}\,=\,Z^{\,n}\,$ for every cell,
given the right-hand side integral $\,I_{\partial\zeta}\,$.
Note that we are actually solving Poisson's equation
$\;{\Delta p}\,=\,{(\,\varrho_{\infty}/\tau)}\,{div\,\vec{u}}\;$
in integral form.
\newline
With the new face pressure gradient computed in following the lines of
section~\ref{S:sec3} the face values of $\,\vec{u}\,$ are updated as
$\,\vec{u}\,-\,(\tau/\varrho_{\infty})\,grad\,p\,$.

After each SOR cycle, continuity of $\,grad\;p\,$
at the cell faces is restored by updating the port values
of $\,p\,$ according to the instructions of the last section.
The loop of processes is reiterated until
$\sum \nolimits_{\zeta}\,I_{\partial\zeta}\,<\,\epsilon\,$
for a suitable bound $\,\epsilon\,$
(\,which happens after a few iterations for appropriate choices\,).
\vspace{-0.2cm}
\section{Convection in coaxial line}\label{S:sec5}
\vspace{-0.2cm}
We illustrate the approach in a stalwart application
by displaying the results of simulations with coaxial line RL100-230
under high power operating conditions
(\,which are realistically inferred from a ion cyclotron resonance
heating \emph{ICRH} experiment in plasma physics\,).

The inner and outer conductors of diameters 100 mm and 230 mm
are made of copper and aluminium, respectively, and the rigid
line is filled with air at atmospheric pressure.
The heating process has been simulated from standby to steady state
CW operation, at frequency 100 MHz and 160 kW transmitted power,
for horizontal position of the line and with outer conductor cooled
at 40 degrees Celsius.

\begin{figure}[!h]
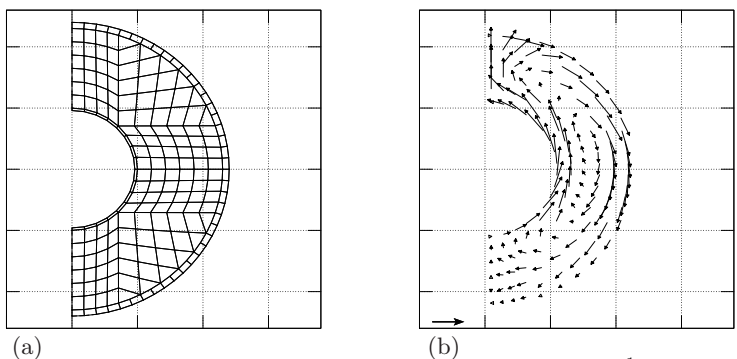
\centering
\vspace{-0.4cm}
\setlength{\unitlength}{1.cm}
\begin{pspicture}(0.0,0.0)(8.0,5.8)\centering
\psset{xunit=0.9cm,yunit=0.9cm}
\rput(0.50,3.70){\includegraphics[scale=0.80,clip=0]{mesh.EPSF}}
\rput(-0.90,1.05){\small (a)}
\rput(6.60,3.70){\includegraphics[scale=0.80,clip=0]{flow.EPSF}}
\rput(5.25,1.05){\small (b)}
\end{pspicture}
\vspace{-1.0cm}
\caption{\textsl{(a) mesh (b) velocity profile
[reference arrow: 0.1$\,m s^{-1}\,$]\qquad\qquad\qquad\hfill}}
\label{F:4}
\vspace{-0.2cm}
\end{figure}

Figure~\ref{F:4}~b displays the computed air flow profile (vertical section)
in steady state, which is attained some minutes after power-on.
Visibly, the natural convection pattern is nicely developed.

Our computations have been carried out with a 3D-mesh of 10 layers
in axial direction, over 200 millimeters of line, the transverse cross
section of which is displayed in figure~\ref{F:4}~a.
At the metallic interfaces no-slip boundary conditions are
implemented and free-slip conditions at all other boundaries.

Simultaneously, and using the same mesh cell system, a Maxwell 
field TLM algorithm was run to provide the heat sources.

\vspace{-0.2cm}
\section{Conclusion}
\vspace{-0.2cm}
A prototypical implementation of the \textsc{Oberbeck-Boussinesq}
approximation to viscous flow has been presented in this paper,
which demonstrates the fundamental fitness of the DSC approach for
fluid dynamic computations.
DSC schemes thus significantly transcend the range of application of
the TLM method from which they descend. \newline
A next natural step in the line of this study is the implementation
of turbulence models which are compatible with the \textsc{Boussinesq}
approach, such as the $\,k-\epsilon\,$ model \cite{ATP}\,, first of all.
Another objective is the incorporation of compressible flow. \newline
-~We hope this paper stimulates some interest into joint further
investigation in these directions.

\hrulefill
\newline
\textsc{Spinner} GmbH. M\"unchen; Aiblinger Str. 30, DE-83620 Westerham
\newline
E-mail address:\; s.hein@spinner.de 
\end{document}